\documentclass[10pt]{article}

\usepackage{amssymb,amsthm,amsmath}
\usepackage[cp1250]{inputenc}
\usepackage{amssymb}

\newcommand{\p}{\mathbb{P}}

\newcommand{\E}{\mathbb{E}}

\newcommand{\R}{\mathbb{R}}
\newcommand{\Z}{\mathbb{Z}}

\newtheorem{theorem}{Theorem}
\newtheorem{lemma}{Lemma}
\newtheorem{cor}{Corollary}

\newtheorem{prop}{Proposition}

\newcommand*{\ind}[1]{\mathbf{1}_{\{#1\}}}

\title{A few remarks on the operator norm of random Toeplitz
matrices}
\author{Rados{\l}aw Adamczak\\ University of Warsaw \thanks{University of Warsaw, Institute of Mathematics,
ul. Banacha 2, 02-097 Warszawa, Poland, e-mail:
radamcz@mimuw.edu.pl.}\;\thanks{This work was completed during the
author's stay at the University of Alberta, Edmonton, Canada
within the PIMS postdoctoral fellowship.}}
\begin{document}
\maketitle
\begin{abstract}
We present some results concerning the almost sure behaviour of
the operator norm or random Toeplitz matrices, including the law
of large numbers for the norm, normalized by its expectation (in
the i.i.d. case). As tools we present some concentration
inequalities for suprema of empirical processes, which are
refinements of recent results by Einmahl and Li.\\


\end{abstract}
\section{Introduction}
For a sequence of independent real random variables
$X_0,X_1,\ldots$, consider the associated sequence of random
Toeplitz matrices, given by
\begin{displaymath}
T_n = \left[
\begin{array}{cccccc}
X_0 & X_1 & X_2 & \ldots &X_{n-2} & X_{n-1}\\
X_1 & X_0 & X_1 &  & &X_{n-2}\\
X_2 & X_1 & X_0 &&&\vdots \\
\vdots &&&\ddots&&\vdots\\
X_{n-2} &&&&X_0&X_1\\
 X_{n-1} & X_{n-2} &\ldots & \ldots & X_1
&X_0
\end{array}
\right] = [X_{|i-j|}]_{1 \le i,j \le n}.
\end{displaymath}

In two recent articles, Bryc, Dembo and Jiang \cite{BDJ} and
independently Hammond and Miller \cite{HM} addressed the questions
raised previously in Bai's paper \cite{B}, concerning the
asymptotic behaviour of the spectral measure of random Toeplitz
matrices generated by i.i.d. mean zero, variance one random
variables. They showed that the spectral measure of $n^{-1/2} T_n$
with probability one converges weakly to a nonrandom measure,
independent of the distribution of the underlying i.i.d. sequence.
Their results were later complemented by Bose and Sen \cite{BS}
and Meckes \cite{M} who analyzed the asymptotic behaviour of the
spectral norm of a random Toeplitz matrix. Bose and Sen proved a
law of large numbers in the case of i.i.d. random variables with
positive mean and finite variance. Meckes showed that if all
$X_i$'s are centered and uniformly subgaussian then $\E\|T_n\|
\simeq \sqrt{n\log n}$ and the same holds with probability one for
$\|T_n\|$, provided that $X_i$'s have some concentration of
measure property.

\paragraph{}
In the present article we are interested in extending the results
of \cite{M} to i.i.d. mean zero random variables with finite
variance as well as to independent non-identically distributed
variables, which satisfy some moment bounds. In particular, we
prove that if $X_0,X_1,X_2,\ldots$ are i.i.d. random variables,
$\E X_0 = 0$ and $0< \E X_0^2 < \infty$, then
\begin{displaymath}
\lim_{n\to \infty}\frac{\|T_n\|}{\E\|T_n\|} = 1\; {\rm a.s.}
\end{displaymath}
This is the statement of Theorem \ref{iid_case} in Section
\ref{section_iid}.

The above result, together with some additional arguments, may be
used to prove that if $X_0,X_1,X_2,\ldots$ are i.i.d. random
variables, then
\begin{align}\label{firstlimsup}
\limsup_{n \to \infty}\frac{\|T_n\|}{\sqrt{n\log n}} < \infty\;
{\rm a.s.},
\end{align}
if and only if $\E X_0 = 0$ and $\E X_0^2 < \infty$, which is the
content of Corollary \ref{main_cor}, Section \ref{section_iid}.

\paragraph{}
In the framework of  independent, but not necessarily identically
distributed variables, we give sufficient conditions for
(\ref{firstlimsup}) to hold (Proposition \ref{non_iid}, Section
\ref{section_non_iid})

\paragraph{}
 Our tools come from the classical
probability in Banach spaces theory and are based on concentration
inequalities for sums of independent Banach space valued
variables. In particular, the proofs are inspired by a recent
paper by Einmahl and Li \cite{EL}, in which a very general law of
the iterated logarithm for sums of i.i.d. variables in Banach
spaces has been characterized. The other crucial ingredient is the
original argument by Meckes which is applied conditionally.

\paragraph{}
The organization of the paper is as follows. In the next section
we present the basic concentration results for sums of independent
random variables. Next, in Section \ref{macierze}, we pass to the
original problem, concerning random matrices. First we obtain some
useful estimates for the expected operator norm and then we prove
results on the almost sure behaviour. We conclude with two
exponential inequalities.

\paragraph{} In the article we will deal both with absolute constants and constants depending on some parameters,
which will be explicitly stated or clear from the context. In both
cases the values of constants may change from line to line.

\paragraph{Acknowledgements} The author would like to thank Prof. Mark
Meckes and Prof. W{\l}odzimierz Bryc for introducing him to the
problem and valuable remarks concerning the early version of the
article.

\section{Concentration of measure for suprema of empirical processes\label{concentration}}

In this section we introduce the basic tool required for the
analysis of the almost sure behaviour of $\|T_n\|$, i.e.
concentration inequalities for sums of independent Banach space
valued random variables. The following exposition is motivated by
the work of Einmahl and Li. In particular, our Theorem
\ref{FukNagaevThm} is a version of Theorem 3.1. from \cite{EL}.

All the concentration results presented below are based on the
celebrated inequality by Talagrand \cite{T}. We will use a recent
version by Klein and Rio, which gives the optimal constant 2 in
front of the variance. We would also like to mention that since
the natural normalization of the operator norm of Toeplitz
matrices is $\sqrt{n\log n}$, which is much bigger than the
classical law of the iterated logarithm normalization, this
optimal value of the constant is not important for our
applications (in fact any other constant would work). However,
since inequalities for sums of Banach space valued variables with
optimal value of this constant are of independent interest and
obtaining them does not require much additional work, we present
the results of this section in a slightly larger generality, than
necessary for subsequent applications.

\paragraph{}Let us first recall the result due to Klein and Rio

\begin{theorem}[\cite{KR}, Theorems 1.1, 1.2]\label{KleinRio}
Let $X_1,X_2,\ldots,X_n$ be independent random variables with
values in a measurable space $(S,\mathcal{S})$ and let
$\mathcal{F}$ be a countable class of measurable functions $f
\colon \mathcal{S}\to [-M,M]$, such that for all $i$, $\E f(X_i)
=0$. Consider the random variable
\begin{displaymath}
Z =\sup_f \sum_{i=1}^n f(X_i).
\end{displaymath}
Then, for all $t \ge 0$,
\begin{displaymath}
\p(Z \ge \E Z +t) \le \exp\Big[-\frac{t^2}{2(\sigma^2 + 2M\E Z) +
3Mt}\Big]
\end{displaymath}
and
\begin{displaymath}
\p(Z \le \E Z - t) \le \exp\Big[-\frac{t^2}{2(\sigma^2 + 2M\E Z) +
3Mt}\Big],
\end{displaymath}
where
\begin{displaymath}
\sigma^2 = \sup_{f\in \mathcal{F}}\sum_{i=1}^n \E f(X_i)^2.
\end{displaymath}
\end{theorem}

The following corollary was derived from the above theorem by
Einmahl and Li \cite{EL} (they stated only the first part,
concerning the upper tail but their proof yields also the other
inequality).

\begin{cor}\label{CLT_type}
Let $X_1,\ldots,X_n$ be independent mean zero random variables
with values in a separable Banach space $B$. Assume that
$\|X_i\|_\infty \le M$ for every $i$ and let
\begin{displaymath}
\sigma^2 = \sup_{f \in B^\ast, \|f\|\le 1} \sum_i \E f(X_i)^2.
\end{displaymath} Consider the random variable
\begin{displaymath}
Z =\|\sum_{i=1}^n X_i\|.
\end{displaymath}
Then for all $0 < \eta \le 1$, $\delta > 0$ there exists a
constant $K = K(\eta,\delta)$, such that for all $t \ge 0$,
\begin{align*}
\p(\max_{1\le k \le n} \|\sum_{i=1}^ k X_i\| \ge (1+\eta)\E Z + t)
\le& \exp\Big[-\frac{t^2}{2(1+\delta)\sigma^2}\Big] \\
&+ \exp\Big[-\frac{t}{K\max_i\|X_i\|_\infty}\Big]
\end{align*}
and
\begin{align*}
\p(Z \le (1-\eta)\E Z - t) \le
\exp\Big[-\frac{t^2}{2(1+\delta)\sigma^2}\Big] +
\exp\Big[-\frac{t}{K\max_i\|X_i\|_\infty}\Big].
\end{align*}
\end{cor}

Using the above result, Einmahl and Li \cite{EL} obtained an
infinite dimensional version of the Fuk-Nagaev inequality for the
deviation of maxima of partial sums of independent random
variables \cite{FN}. For our purposes we will need a slightly
refined version of their estimate, given in the following

\begin{theorem}\label{FukNagaevThm}
In the setting of Corollary \ref{CLT_type} assume that for some $p
\ge 1$ and all $i$, $\E\|X_i\|^p < \infty$. Then for all $0 < \eta
\le 1$, $\delta > 0$ there exists a constant $C =
C(p,\eta,\delta)$, such that for all $t \ge 0$,
\begin{align}\label{FukNagaev1}
\p&\Big(\max_{1\le k \le n} \|\sum_{i=1}^k X_i\| \ge (1+
\eta)\E Z + t\Big)\nonumber\\
&\le\exp\Big[-\frac{t^2}{2(1+\delta)\sigma^2}\Big]+ C \E\max_{1\le
i \le n}\|X_i\|^p/t^p
\end{align}
and
\begin{align}\label{FukNagaev2}
\p&\Big( Z \le (1-
\eta)\E Z - t\Big)\nonumber\\
&\le\exp\Big[-\frac{t^2}{2(1+\delta)\sigma^2}\Big]+ C\E \max_{1\le
i\le n} \|X_i\|^p/t^p.
\end{align}
\end{theorem}

The proof we give below follows quite closely the original
argument of Einmahl and Li and is based on the
Hoffmann-J{\o}rgensen inequality. The difference is that we use it
in the moment version, which in our opinion simplifies the
argument and also allows to replace $\sum_i \E\|X_i\|^p$,
appearing in the original estimate by Einmahl and Li, by $\E\max_i
\|X_i\|^p$ (this will be useful when dealing with Toeplitz
matrices generated by non i.i.d. variables). A somewhat similar
application to moment estimates appeared in \cite{GLZ}, whereas
\cite{A} contains an analogous derivation in the context of
exponential inequalities for unbounded empirical processes (its
final version was also influenced by results of Einmahl and Li).
We state it at the end of the present section.

\begin{proof}
 For $\varepsilon = \varepsilon(\delta)
> 0$ (whose value will be determined later) consider truncated
variables $Y_i = X_i\ind{\|X_i\|\le \rho}$. Set also $R_i = X_i -
Y_i$.

We have
\begin{align*}
\max_{1\le k \le n}\|\sum_{i=1}^k X_i\| \le \max_{1\le k \le
n}\|\sum_{i=1}^k (Y_i - \E Y_i)\| + \max_{1\le k \le
n}\|\sum_{i=1}^k (R_i - \E R_i)\|
\end{align*}
and
\begin{align*}
\|\sum_{i=1}^n X_i\| \ge \|\sum_{i=1}^n (Y_i - \E Y_i)\| -
\|\sum_{i=1}^n (R_i- \E R_i)\|,
\end{align*}
where we used the fact that $\E Y_i + \E R_i = 0$.

Similarly, using Jensen's inequality, we obtain
\begin{align*}
\E \|\sum_{i=1}^n (Y_i- \E Y_i)\| - 2\E \|\sum_{i=1}^n R_i\| \le \E
Z \le \E \|\sum_{i=1}^n (Y_i- \E Y_i)\| + 2\E \|\sum_{i=1}^n R_i\|.
\end{align*}
Thus, setting $A = \E\max_{1\le k \le n}\|\sum_{i=1}^k R_i\|$, we
get
\begin{align}\label{aux1}
\p(&\max_{1\le k \le n}\|\sum_{i=1}^k X_i\| \ge (1 + \eta)\E Z +
t)\nonumber \\
 \le& \p(\max_{1\le k\le n} \|\sum_{i=1}^k (Y_i  - \E Y_i)\| \ge (1
+ \eta)\E Z + (1-\varepsilon)t) \nonumber \\
&+ \p(\max_{1\le k\le n}\|\sum_{i=1}^k (R_i - \E R_i)\| \ge
\varepsilon
t)\nonumber\\
 \le& \p(\max_{1\le k\le n} \|\sum_{i=1}^k (Y_i  - \E Y_i)\| \ge (1
+ \eta)\E\|\sum_{i=1}^k (Y_i  - \E Y_i)\| - 4 A + (1-\varepsilon)t) \nonumber \\
&+ \p(\max_{1\le k\le n}\|\sum_{i=1}^k (R_i - \E R_i)\| \ge
\varepsilon t)
\end{align}
and similarly
\begin{align}\label{aux2}
\p(&Z \le (1-\eta)\E Z - t)  \nonumber\\
\le & \p(\|\sum_{i=1}^n (Y_i - \E Y_i)\|  \le (1- \eta)\E Z -
(1-\varepsilon)t )\nonumber \\
&+ \p(\|\sum_{i=1}^n (R_i - \E R_i)\| \ge \varepsilon t) \nonumber\\
\le& \p(\|\sum_{i=1}^n (Y_i - \E Y_i)\|  \le (1-
\eta)\E\|\sum_{i=1}^n (Y_i - \E Y_i)\| -
(1-\varepsilon)t + 2 A)\nonumber \\
&+ \p(\|\sum_{i=1}^n (R_i - \E R_i)\| \ge \varepsilon t).
\end{align}
We would like to apply Corollary \ref{CLT_type} to bound first
summands on right hand sides of the above inequalities and the
Markov inequality to bound the deviation of the remainder. Thus
the aim is now to find a good truncation level $\rho$, which would
ensure that $A \simeq \varepsilon t$ and both
\begin{displaymath}
e^{-t(1-\mathcal{O}(\varepsilon))/K\max_i\|Y_i-\E Y_i\|_\infty},
\p(\max_{1\le k\le n} \|\sum_{i=1}^k (R_i - \E R_i)\| \ge
\varepsilon t)
\end{displaymath}
are of order $t^{-p}\E\max_{i\le n}\|X_i\|^p$, where $K =
K(\eta,\delta)$ is the constant from Corollary \ref{CLT_type}.

Recall now the classical Hoffmann-J{\o}rgensen inequality (e.g.
\cite{LT}, Chapter 6, Proposition 6.8., inequality (6.8)), which
asserts that
\begin{displaymath}
\E\max_{1\le k\le n}\|\sum_{i=1}^k R_i\|^p \le 2\cdot
4^p\E\max_{1\le i\le n} \|R_i\|^p + 2(4t_0)^p,
\end{displaymath}
where
\begin{displaymath}
t_0 = \inf\{t >0 \colon \p(\max_{1\le k\le n}\|\sum_{i=1}^k R_i\|
> t) \le (2\cdot 4^p)^{-1}\}.
\end{displaymath}
The idea (following \cite{GLZ}) is to truncate at a level which
would allow to take $t_0 = 0$. Let us thus set
\begin{align}\label{choice_of_rho}
\rho^p = 2\cdot 4^p\cdot\E\max_{1\le i \le n} \|X_i\|^p.
\end{align}
Then $\p(\max_{1\le k\le n}\|\sum_{i=1}^k R_i\| > 0) \le
\p(\max_{1\le i \le n} \|X_i\| > \rho) \le (2\cdot 4^p)^{-1}$ by
the definition of $R_i$ and the Markov inequality. Thus
\begin{displaymath}
\E\max_{1\le k \le n} \|\sum_{i=1}^k R_i\|^p \le 2\cdot4^p
\cdot\E\max_{1\le i\le n} \|X_i\|^p
\end{displaymath}
and also
\begin{displaymath}
\E\max_{1\le k \le n} \|\sum_{i=1}^k R_i\| \le 2^{1/p}\cdot 4
\cdot (\E\max_{1\le i\le n} \|X_i\|^p)^{1/p}.
\end{displaymath}

Notice, that we can assume that $\E\max_{1\le i\le n}\|X_i\|^p\le
\varepsilon^p t^p$, since otherwise the statement of the corollary
becomes trivial ($\varepsilon$ is a function of $\delta$, so it
can be used to define $C$ such that right hand sides of the
inequalities in question are greater than 1, if $\E\max_{1\le i\le
n}\|X_i\|^p> \varepsilon^p t^p$).

Therefore we have
\begin{align}\label{pomocniczy}
A = \E\max_{1\le k \le n} \|\sum_{i=1}^k R_i\| \le 8t\varepsilon.
\end{align}
Also, by Jensen's inequality
\begin{displaymath}
\E\max_{1\le k\le n}\|\sum_{i=1}^k (R_i - \E R_i)\|^p \le
2^p\E\max_{1\le k\le n}\|\sum_{i=1}^k R_i\|^p \le 2\cdot 8^p
\E\max_{1\le i \le n} \|X_i\|^p,
\end{displaymath}
so Markov's inequality gives
\begin{align*}
\p(\max_{1\le k\le n}\|\sum_{i=1}^k (R_i - \E R_i)\| \ge \varepsilon
t) \le 2\cdot 8^p\varepsilon^{-p}\frac{\E\max_{1\le i \le
n}\|X_i\|^p}{t^p}.
\end{align*}

\noindent Combining the above estimates with inequalities
(\ref{aux1}), (\ref{aux2}) and (\ref{pomocniczy}) gives
\begin{align*}
\p(&\max_{1\le k \le n}\|\sum_{i=1}^k X_i\| \ge (1 + \eta)\E Z +
t)\nonumber\\
\le&\p(\max_{1\le k\le n} \|\sum_{i=1}^k (Y_i  - \E Y_i)\| \ge (1
+ \eta)\E\|\sum_{i=1}^k (Y_i  - \E Y_i)\|  + (1-33\varepsilon)t) \nonumber \\
&+2\cdot 8^p\varepsilon^{-p}\frac{\E\max_{1\le i \le
n}\|X_i\|^p}{t^p}
\end{align*}
and
\begin{align*}
\p(&Z \le (1-\eta)\E Z - t)  \nonumber\\
\le& \p(\|\sum_{i=1}^n (Y_i - \E Y_i)\|  \le (1-
\eta)\E\|\sum_{i=1}^n (Y_i - \E Y_i)\| -
(1-17\varepsilon)t)\nonumber \\
&+2\cdot 8^p\varepsilon^{-p}\frac{\E\max_{1\le i \le
n}\|X_i\|^p}{t^p}.
\end{align*}

Let us notice that for every $f$ in the unit ball of $B^\ast$ we
have $\E [f(Y_i - \E Y_i)]^2 = \E [f(Y_i) - \E f(Y_i)]^2 \le \E
f(Y_i)^2 \le \E f(X_i)^2$. Thus, by Corollary \ref{CLT_type} we can write

\begin{align*}
\p(&\max_{1\le k \le n}\|\sum_{i=1}^k X_i\| \ge (1 + \eta)\E Z +
t), \quad \p(Z \le (1-\eta)\E Z - t) \nonumber\\
\le &
\exp\Big[-\frac{t^2(1-33\varepsilon)^2}{2(1+\delta)\sigma^2}\Big]
+ \exp\Big[-\frac{(1-33\varepsilon)t}{2K(\eta,\delta) \rho}\Big] +
2\cdot 8^p\varepsilon^{-p}\frac{\E\max_{1\le i \le
n}\|X_i\|^p}{t^p}.
\end{align*}

Since $\exp(-x) \le C_p x ^{-p}$, we have for $\varepsilon <
1/66$,
\begin{align*}
\exp\Big[-\frac{(1- 33\varepsilon)t}{2K(\eta,\delta)\rho}\Big]
\le& C_p [K(\eta,\delta)\rho]^pt^{-p} = C_p K(\eta,\delta)^p\cdot
2\cdot 4^p \frac{\E\max_{1\le i \le n}\|X_i\|^p}{t^p}\\
\le& C(p,\eta,\delta)\frac{\E\max_{1\le i \le n}\|X_i\|^p}{t^p}.
\end{align*}

To finish the proof it is enough to apply the last two
inequalities with $\delta/2$ instead of $\delta$ and $\varepsilon$
small enough to ensure that
\begin{displaymath}
(1 + \delta/2)(1 - 33\varepsilon)^{-2} \le (1 + \delta).
\end{displaymath}
\end{proof}

As announced, we would  also like to state a similar result,
obtained in \cite{A}, valid for variables satisfying stronger,
exponential integrability conditions. We will use it to obtain
concentration inequalities for $\|T_n\|$ (see Proposition
\ref{conc_ineq} in Section \ref{last_section} below).

\begin{theorem}\label{unbounded}
In the setting of Corollary \ref{CLT_type}, assume that for some
$0< \alpha\le 1$ and all $i$, $\|X_i\|_{\psi_\alpha} < \infty$.
For all $0< \eta < 1$ and $\delta
> 0$, there exists a constant $C =C(\alpha,\eta,\delta)$, such
that for all $t \ge 0$,
\begin{align*}
\p(Z \ge (1+\eta)&\E Z + t) \\
&\le \exp\Big[-\frac{t^2}{2(1+\delta)\sigma^2}\Big] +
3\exp\Big[-\Big(\frac{t}{C\|\max_{i}\|X_i\|\|_{\psi_\alpha}}\Big)^{\alpha}\Big]
\end{align*}
and
\begin{align*}
\p(Z \le (1-\eta)&\E Z - t) \\
&\le \exp\Big[-\frac{t^2}{2(1+\delta)\sigma^2}\Big] +
3\exp\Big[-\Big(\frac{t}{C\|\max_{i}\|X_i\|\|_{\psi_\alpha}}\Big)^{\alpha}\Big].
\end{align*}
\end{theorem}

\section{Random Toeplitz matrices\label{macierze}}

In this part we will focus on the initial problem and prove
several results, concerning the operator norm of random Toeplitz
matrices. The consecutive sections will be devoted to properties
of the expectation, almost sure behaviour and concentration
inequalities.

\subsection{Expectation estimates}
We will start with some estimates of the expected operator norm.
In the proof of our results we use an argument introduced by
Meckes, but apply it conditionally.

Following \cite{M}, let us first note that $T_n$ is a submatrix of
the infinite Laurent matrix
\begin{displaymath}
L_n = [X_{i-1}\ind{|i-j|\le n-1}]_{i,j\in \Z},
\end{displaymath}
so $\|T_n\| \le \|L_n\|$, where by $\|L_n\|$ we denote the
operator norm of $L_n$ acting in the standard way on $\ell_2(\Z)$.
If we use the Fourier basis to identify $\ell_2(\Z)$ with
$L_2[0,1]$, it turns out that $L_n$ corresponds to a
multiplication operator, with the multiplier
\begin{displaymath}
f(t) = X_0 + 2\sum_{j=1}^{n-1} \cos(2\pi j t)X_j
\end{displaymath}
and thus $\|L_n\| = \|f\|_\infty$.

\begin{theorem}\label{expectation_from_above}
If $X_0,X_1,\ldots$ are independent mean zero random variables,
then
\begin{displaymath}
\E \|T_n\| \le \E\|L_n\|\le C \sqrt{\sum_{i=0}^{n-1} \E
X^2_i}\sqrt{\log n},
\end{displaymath}
where $C$ is an absolute constant.
\end{theorem}

\begin{proof}
By the arguments preceding the formulation of the theorem and the
standard symmetrization inequality, we have
\begin{displaymath} \E \|T_n\| \le \E\|L_n\| \le
2\E\sup_{t\in[0,1]} |Y_t|,
\end{displaymath}
where
\begin{displaymath}
Y_t = X_0\varepsilon_0 + 2\sum_{j=1}^{n-1} \cos(2\pi j
t)\varepsilon_jX_j
\end{displaymath}
and $\varepsilon_1,\varepsilon_2,\ldots$ is a sequence of i.i.d.
Rademacher variables, independent of $X_j$'s. Using the fact that
Rademacher variables are subgaussian and Dudley's entropy bound,
conditionally on $X_i$'s, one gets
\begin{align}\label{conditional}
\E_\varepsilon \sup_{t\in [0,1]} |Y_t| \le |Y_0| + C\int_0^\infty
\sqrt{\log N([0,1],d,s)}ds,
\end{align}
where $C$ is an absolute constant,
\begin{displaymath}
d(t_1,t_2) = \|Y_{t_1} - Y_{t_2}\|_2 = 2\sqrt{\sum_{j=1}^{n-1}
X_j^2[\cos(2\pi j t_1) - \cos(2\pi j t_2)]^2 }
\end{displaymath}
and $N([0,1],d,s)$ is the smallest cardinality of an $s$-net in
$[0,1]$ with the metric $d$ (see e.g. \cite{VW}, chapter 2.2).

 We have ${\rm diam}([0,1],d) \le 4\sqrt{\sum_{j=1}^{n-1}
X_j^2} =: D$. Moreover
\begin{displaymath}
d(t_1,t_2) \le 4\pi\sqrt{\sum_{j=1}^{n-1} j^2 X_j^2}|t_1-t_2|,
\end{displaymath}
so
\begin{displaymath}
N([0,1],d,s) \le C\frac{A}{s},
\end{displaymath}
for $s \le D$, where $A = \sqrt{\sum_{j=1}^{n-1} X_j^2j^2}$.

Therefore
\begin{align*}
\int_0^\infty &\sqrt{\log N([0,1],d,s)}ds \le \int_0^D\sqrt{\log
\frac{CA}{s}}ds\\
&= \frac{CA}{\sqrt{2}}\int_{\sqrt{2\log{(CA/D)}}}^\infty t^2
e^{-t^2/2}dt\\
&\le \frac{CA}{\sqrt{2}}\Big(\sqrt{2\log{(CA/D)}} +
\sqrt{2\pi}\Big)e^{-\log{(CA/D)}}\\
&\le D\sqrt{\log{(CA/D)}} +\sqrt{\pi}D.
\end{align*}

We have $A/D \le \sqrt{\sum_{j=1}^{n-1} j^2} \le Cn^{3/2}$, so the
above inequality  together with (\ref{conditional}) implies
\begin{displaymath}
\E_\varepsilon \sup_{t\in [0,1]} |Y_t| \le |Y_0| +
C\sqrt{\sum_{j=1}^{n-1} X_j^2}\sqrt{\log n}.
\end{displaymath}
Taking into account that $Y_0 = X_0\varepsilon_0 +
2\sum_{j=1}^{n-1}X_j\varepsilon_j$, integrating with respect to
the variables $X_j$'s and using Jensen's inequality, we get
\begin{displaymath}
\E\|L_n\|\le C\sqrt{\sum_{j=0}^{n-1} \E X_j^2}\sqrt{\log n}.
\end{displaymath}
\end{proof}

The following lemma is a partial converse to the above theorem. It
also relies on arguments from \cite{M}.
\begin{lemma}\label{necessity}
If $X_0, X_1,\ldots,$ are i.i.d. and
\begin{displaymath}
\limsup_{n\to \infty} \frac{\E\|T_n\|}{\sqrt{n\log n}} < \infty,
\end{displaymath}
then $\E X_i^2 <\infty$.
\end{lemma}

\begin{proof}
Without loss of generality we can assume that $X_i$'s are
symmetric with positive variance. Thus for $R$ large enough
$X_j\ind{|X_j|\le R}$ is not identically zero.

Let $\varepsilon_1,\varepsilon_2,\ldots$ be i.i.d. Rademacher
variables, independent of $X_i$'s. By assumption, for some $C <
\infty$,
\begin{displaymath}
\E\|T_n\| \le C\sqrt{n\log n}.
\end{displaymath}
By another argument from \cite{M} (see the proof of Theorem 3
therein), we have
\begin{displaymath}
\|T_n\| \ge \sup_{t\in [0,1]}\Big|X_0 +\sum_{j=1}^{n-1}
2(1-j/n)X_j\cos(2\pi jt)\Big|.
\end{displaymath}
Thus, using symmetry of $X_j$'s and the contraction principle, we
obtain that for each $R < \infty$,
\begin{displaymath}
\E\|T_n\| \ge \E_X\E_\varepsilon\sup_{t\in[0,1]}
\Big|\sum_{j=0}^{\lfloor n/2\rfloor} \varepsilon_j
X_j\ind{|X_j|\le R}\cos(2\pi jt)\Big|.
\end{displaymath}
Now, by Proposition 6 of \cite{M}, we obtain for $R$ large enough
\begin{displaymath}
\E\|T_n\| \ge c\E\|(a_j)_{j=1,\ldots,\lfloor
n/2\rfloor}\|_2\sqrt{\log\frac{\|(a_j)_{j=1,\ldots,\lfloor
n/2\rfloor}\|_2}{\|(a_j)_{j=1,\ldots,\lfloor n/2\rfloor}\|_4}},
\end{displaymath}
where $a_j = X_j \ind{|X_j|\le R}$.  Let $A = \E
X_0^2\ind{|X_0|\le R}$. By the law of large numbers, for $n$ large
enough (depending on $R$), the random variable on the right hand
side is with probability at least $1/2$ greater than or equal to
\begin{displaymath}
\sqrt{n}\sqrt{A/2}\sqrt{\log{\frac{n^{1/2}\sqrt{A/2}}{n^{1/4}R}}}\ge
\sqrt{n}\sqrt{A/2} \sqrt{\log{(n^{1/8}\sqrt{A/2})}}.
\end{displaymath}
Thus, for $n$ large enough (depending on $R$ and $A$), we have
\begin{displaymath}
\sqrt{nA} \sqrt{\log n^{1/16}} \le \sqrt{nA}
\sqrt{\log{(n^{1/8}\sqrt{A})}} \le \sqrt{2}Cc^{-1} \sqrt{n\log n},
\end{displaymath}
which gives $A \le 32 C^2c^{-2}$. Since $R$ is arbitrary, it
implies that $\E X_0^2 \le 32C^2c^{-2}$.
\end{proof}

\subsection{Almost sure behaviour\label{asbehaviour}}

Before we proceed to the analysis of the almost sure behaviour of
the operator norm of $T_n$, let us introduce some notation that
will be used throughout the rest of the paper.

For every positive integer $n$ and $i=0,\ldots,n-1$, let
$A_i^{(n)} = [\ind{|j -l|=i}]_{j,l \le n}$. We thus have
\begin{align}\label{sum_representation}
T_n = \sum_{i=0}^{n-1} X_i A_i^{(n)}.
\end{align}
This basic relation will allow us to apply to $T_n$ the results
for sums of independent random variables, presented in Section
\ref{concentration}.

\subsubsection{The non i.i.d. case \label{section_non_iid}}

As the first example of application of Theorem \ref{FukNagaevThm}
we will prove the following

\begin{prop}\label{non_iid}
If $X_1,X_2,\ldots$ are independent (not necessarily i.i.d.), $\E
X_n = 0$ and for some $\alpha > 1$, $\sup_n \E |X_n|^2[\log\log
(|X_n|\vee e^e)]^\alpha  < \infty$, then with probability 1,
\begin{displaymath}
\limsup_{n\to \infty}\frac{\|T_n\|}{\sqrt{n\log n}} =
\limsup_{n\to \infty} \frac{\E\|T_n\|}{\sqrt{n\log n}} < \infty\;
\end{displaymath}
\end{prop}

\begin{proof}
The inequality follows immediately from Theorem
\ref{expectation_from_above}. We will first prove, that
\begin{align}\label{real_version}
\limsup_{n\to \infty}\frac{\|T_n\|}{\sqrt{n\log n}} \le
\limsup_{n\to \infty} \frac{\E\|T_n\|}{\sqrt{n\log n}} \;{\rm a.s.}
\end{align}
It turns out that to obtain the inverse of (\ref{real_version}),
we will need a corresponding result for Laurent matrices, i.e.
\begin{align}\label{real_version_L}
\limsup_{n\to \infty}\frac{\|L_n\|}{\sqrt{n\log n}} \le
\limsup_{n\to \infty} \frac{\E\|L_n\|}{\sqrt{n\log n}} \;{\rm
a.s.}
\end{align}
The proofs of both statements are very similar.  Therefore we will
present in detail only the argument for $T_n$ and then we will
indicate all the changes one has to make in order to obtain
(\ref{real_version_L}).

\paragraph{}
Let us denote $A = \limsup \E\|T_n\|/\sqrt{n\log n}$.
 To obtain (\ref{real_version}), it is enough to show that for all $\varepsilon >
0$,
\begin{align}\label{prob_version}
\lim_{n\to \infty} \p\Big(\sup_{k>n} \|T_k\|/\sqrt{k\log k} > A +
\varepsilon\Big) = 0
\end{align}

We have for $\theta > 1$,
\begin{align}\label{Toeplitz_special}
\p\Big(\sup_{k > n} &\|T_k\|/\sqrt{k\log k} > A +
\varepsilon\Big) \nonumber\\
&\le\sum_{k \ge \lfloor \log n /\log \theta \rfloor}
\p(\max_{\theta^k \le i < \theta^{k+1}} \|T_i\| > (A +
\varepsilon)\sqrt{\theta^k\log \theta^k}\Big)\nonumber \\
&=\sum_{k \ge \lfloor \log n /\log \theta \rfloor}
\p(\|T_{\theta^{k+1}}\| > (A + \varepsilon)\sqrt{\theta^k\log
\theta^k}\Big),
\end{align}
where in the last inequality we used the fact that $\|T_i\|$ is
increasing with $i$ (since $T_{i}$ is a submatrix of $T_{i+1}$).
Now, for $\eta$ small enough and $k$ large enough, we have
\begin{align*}
(A + \varepsilon/2)\sqrt{\theta^k\log \theta^k} &\ge (A +
\varepsilon/4)(1 + 2\eta) \sqrt{\theta^k\log\theta^k}\\
& \ge (1+ 2\eta)
\E\|T_{\theta^{k+1}}\|\sqrt{\theta^{-1}\frac{k\log\theta}{(k+1)\log\theta}},
\end{align*}
which for $\theta$ close to one and $k$ large enough is greater
than or equal to $(1+\eta)\E\|T_{\theta^{k+1}}\|$. Thus, to prove
(\ref{prob_version}), it is enough to show that
\begin{displaymath}
\sum_{k=1}^\infty\p\Big( \|T_{\theta^{k+1}}\| \ge
\E\|T_{\theta^{k+1}}\|(1 + \eta) + \varepsilon \sqrt{\theta^k\log
\theta^k}/2\Big) < \infty
\end{displaymath}
for arbitrary $\theta > 1$ and  $\varepsilon, \eta > 0$. Let us
now denote $\Psi(x) = x[\log\log(x\vee M)]^\alpha$ and $\Phi(x) =
x/[\log\log (x\vee M)]^\alpha$, where $M$ is a parameter and note
that for $M$ large enough, $\Psi$ is convex and nondecreasing and
$\Phi$ is nondecreasing. Moreover for some constant $K$ (depending
on $M$ and $\alpha$) and all $x \ge 0$,
\begin{displaymath}
K+ \Phi(\Psi(x)) \ge x/K.
\end{displaymath}
We will now take into account the relation
(\ref{sum_representation}) to treat $T_n$ as a sum of independent
variables in the space of matrices and apply Theorem
\ref{FukNagaevThm}. Let us notice that since $\|A_i^{(n)}\| \le
2$, $\sigma^2$ from this theorem is bounded by
$4\theta^{k+1}\sup_{i<\theta^{k+1}}\E X_i^2 \le K_1\theta^{k+1}$ for
some constant $K_1$ (depending on the supremum in the assumption of the proposition). Therefore, applying Theorem \ref{FukNagaevThm}
with $p = 2$, we obtain
\begin{align*}
\p\Big( \|T_{\theta^{k+1}}\| &\ge \E\|T_{\theta^{k+1}}\|(1 + \eta)
+ \varepsilon \sqrt{\theta^k\log \theta^k}/2\Big)\\
& \le e^{-\frac{k\varepsilon^2\log\theta}{K_2}} + K_2\frac{\E\max_{i<
\theta^{k+1}}|X_i|^2}{\varepsilon^2 \theta^k k \log\theta}.
\end{align*}

The series corresponding to the first term on the right hand side
is clearly convergent. As for the second term, by the
integrability assumption on $X_i$'s, we have for some constants
$C_1,C_2$,
\begin{align*}
K^{-1}\E\max_{i< \theta^{k+1}}|X_i|^2 &\le K +
\Phi(\Psi(\E\max_{i<
\theta^{k+1}}|X_i|^2)) \le K + \Phi(\E\sum_{i< \theta^{k+1}}\Psi(|X_i|^2))\\
&\le K + \Phi(\theta^{k+1} C_1) \le K +
\frac{\theta^{k+1}C_2}{\log^\alpha k},
\end{align*}
which shows that also the other series is convergent, as
\begin{displaymath}
\sum_{k=1}^\infty \frac{1}{k\log^\alpha k} < \infty
\end{displaymath}
for $\alpha >1$. This proves (\ref{prob_version}) and thus also (\ref{real_version}).

\paragraph{}
Let us now indicate the changes one has to make in the above
argument, to prove (\ref{real_version_L}). First, let us notice
that the matrices $L_n$ are sums of independent variables with
values in $\mathcal{B}(\ell_2(\Z))$ - the space of all bounded
operators on the Hilbert space $\ell_2(\Z)$. Indeed, we have
\begin{align}\label{Laurent_representation}
L_n = \sum_{i=0}^{n-1}X_i A_i,
\end{align}
where $A_i =[\ind{|i-j| = i}]_{i,j\in \Z}$.

The space $\mathcal{B}(\ell_2(\Z))$ is not separable, but one can
still apply Theorem \ref{KleinRio}, since the norm on this space
is expressible as a supremum over a countable set of functionals,
so one can approximate it by norms in separable spaces. The second
difference is the lack of monotonicity of $\|L_n\|$ with respect
to $n$, which prevents us from using a counterpart of the equality
(\ref{Toeplitz_special}). This is however not a real problem,
since, denoting $A = \limsup \|L_n\|/\sqrt{n\log n}$, one can write
\begin{align*}
\p\Big(\sup_{k > n} &\|L_k\|/\sqrt{k\log k} > A +
\varepsilon\Big) \nonumber\\
&\le\sum_{k \ge \lfloor \log n /\log \theta \rfloor}
\p(\max_{\theta^k \le i < \theta^{k+1}} \|L_i\| > (A +
\varepsilon)\sqrt{\theta^k\log \theta^k}\Big)\nonumber \\
&\le\sum_{k \ge \lfloor \log n /\log \theta \rfloor} \p(\max_{i\le
\theta^{k+1}} \|L_{i}\| > (A + \varepsilon)\sqrt{\theta^k\log
\theta^k}\Big)
\end{align*}
and use the full strength of Theorem \ref{KleinRio}, which gives a
bound for the tail of maxima of partial sums. The rest of the
proof is the same as for $T_n$ (note that we still have $\|A_i\|
\le 2$).

\paragraph{}
What remains to be proved is the inverse inequality, i.e.
\begin{align}\label{opposite}
\limsup_{n\to \infty}\frac{\|T_n\|}{\sqrt{n\log n}} \ge
\limsup_{n\to \infty} \frac{\E\|T_n\|}{\sqrt{n\log n}}.
\end{align}

This can be established again by arguments similar as in
\cite{EL}, Section 4.3 (see also \cite{AKL}, Theorem 7). First,
let us define $X_n' = X_n\ind{|X_n| \le \sqrt{n\log n}}$ and let
$T_n'$ (resp. $L_n'$) be the Toeplitz (resp. Laurent) matrix
generated by $X_n'$. We have
\begin{align}\label{first_limsups_eq}
\limsup_{n\to \infty}\frac{\|T_n\|}{\sqrt{n\log n}} =
\limsup_{n\to \infty} \frac{\|T_n'\|}{\sqrt{n\log n}}\; {\rm a.s.}
\end{align}
and
\begin{align}\label{second_limsups_eq}
\limsup_{n\to \infty}\frac{\E\|T_n\|}{\sqrt{n\log n}} =
\limsup_{n\to \infty} \frac{\E\|T_n'\|}{\sqrt{n\log n}}.
\end{align}

 To prove the first equality, it is enough to notice
that by the Chebyshev inequality
\begin{align*}
\sum_{n} \p(X_n \neq X_n') &\le \sum_{n} \p( |X_n| > \sqrt{n\log n}) \\
&\le \sum_{n} C\frac{\E|X_n|^2[\log\log(|X_n|\vee
e^e)]^\alpha}{n[\log n][\log\log n]^\alpha} < \infty,
\end{align*}
which by the Borel-Cantelli Lemma implies that with probability 1,
$X_n = X_n'$ for large $n$. The equality (\ref{first_limsups_eq})
follows by (\ref{sum_representation}), an analogous formula for
$T_n'$ and the fact that $\|A_i^{(n)}\| \le 2$.

Let us now turn to the proof of (\ref{second_limsups_eq}). Notice
that
\begin{align*}
\Big|\E\|T_n\| - \E\|T_n'\|\Big| &\le \sum_{i=0}^{n-1} \|A_i^{(n)}\|\E|X_i|\ind{|X_i|> \sqrt{i\log i}} \\
&\le 2\sum_{i=0}^n\E|X_i|\ind{|X_i| > \sqrt{i\log i}} \le K\sum_{i=0}^{n-1}\frac{\E |X_i|^2}{\sqrt{i\log i}}\\
&= o(\sqrt{n\log n}),
\end{align*}
which already implies (\ref{second_limsups_eq}).

To prove (\ref{opposite}), it is thus sufficient to show that
$\sup_n \|L_n'\|/\sqrt{n\log n}$ is integrable, since then, by
Fatou's lemma (recall that $\|T_n'\|\le \|L_n'\|$), we get
\begin{displaymath}
\E \limsup_{n} \frac{\|T_n'\|}{\sqrt{n\log n}} \ge \limsup_n
\frac{\E\|T_n'\|}{\sqrt{n\log n}},
\end{displaymath}
which allows us to finish the proof, since by (\ref{real_version}), (\ref{first_limsups_eq}) and (\ref{second_limsups_eq}), we
have
\begin{displaymath}
\limsup_{n} \frac{\|T_n'\|}{\sqrt{n\log n}} \le \limsup_n
\frac{\E\|T_n'\|}{\sqrt{n\log n}}\;{\rm a.s.}
\end{displaymath}

To prove the desired integrability, notice that by Theorem \ref{expectation_from_above},
(\ref{real_version_L}), (\ref{Laurent_representation}), the
inequality $\|A_i\|\le 2$ and the fact that with probability 1,
$X_n = X_n'$ for large $n$, we have
\begin{displaymath}
\sup_{n}\frac{\|L_n'\|}{\sqrt{n\log n}} < \infty \;{\rm a.s.}
\end{displaymath}
Moreover $\E\sup_{n}\|A_n X_n'\|/\sqrt{n\log n} \le 2$. By
Proposition 6.12. from \cite{LT}, this implies that also
\begin{displaymath}
\E\sup_n \frac{\|L_n'\|}{\sqrt{n\log n}} < \infty,
\end{displaymath}
which ends the proof.
\end{proof}

\paragraph{}
We would also like to remark that if one considers higher iterates
of logarithm, the assumptions of Proposition \ref{non_iid} may be
weakened, for instance to $\sup_n \E |X_n|^2[\log\log (|X_n|\vee
M)][\log\log\log (|X_n|\vee M)]^\alpha < \infty$ for some
$\alpha>1$ (the proof is essentially the same, the only difference
is in the involved series). The proposition is however no longer
true if we assume only $\sup_n \E |X_n|^2[\log\log (|X_n|\vee e)]
< \infty$. To see this, notice that since $\|T_n\| \ge |X_{n-1}|$,
by the Borel-Cantelli Lemma a necessary condition for the almost
sure finiteness of the $\limsup$ in question is
\begin{displaymath}
\sum_n \p(|X_n| \ge C\sqrt{n\log n}) < \infty
\end{displaymath}
for some $C < \infty$, which is not true for instance if $X_n =
\pm \sqrt{n[\log n][\log\log\log n]}$ with probability $p_n =
(n[\log n][\log\log n][\log\log\log n])^{-1}$ and $0$ with
probability $1 - 2p_n$.

\subsubsection{The i.i.d. case. Proof of main results. \label{section_iid}}

We will now prove the main results of the article, announced in the
Introduction. We state them again below as Theorem \ref{iid_case}
and Corollary \ref{main_cor}. The proofs go along the same lines
as in the non i.i.d. case. The main tool is Theorem \ref{KleinRio}
(now we can actually use the original formulation from \cite{EL},
with sums of $p$-th moments instead of the $p$-th moment of the
maximum of $|X_i|$'s).

\begin{theorem}\label{iid_case}
Assume that $X_0,X_1,X_2,\ldots$ are i.i.d. random variables. If
$\E X_0 = 0$ and $0< \E X_0^2 < \infty$, then
\begin{displaymath}
\lim_{n\to \infty}\frac{\|T_n\|}{\E\|T_n\|} = 1\; {\rm a.s.}
\end{displaymath}
\end{theorem}

\begin{proof} Just as the proof of Proposition \ref{non_iid}, the
argument is influenced by the work of Einmahl and Li. The first
step is the truncation.

Let us define $X'_n = X_n\ind{|X_n|\le \sqrt{n\log n}}$. We have
\begin{align*}
\sum_n \p(X_n \neq X'_n) & = \sum_n \p(X_0 > \sqrt{n\log
n})\\
&= \E\sum_n \ind{n\log n < X_0^2} \le K \E X_0^2 < \infty
\end{align*}
so with probability 1, $X_n = X'_n$ for large $n$.

By Theorem \ref{expectation_from_above}, we have $\E\|T_n\| \le
C\sqrt{n\log n}$ for some $C <\infty$ (depending on the
distribution of $X_0$). Also, by Theorem 3 of \cite{M},
symmetrization inequalities and Jensen's inequality, we have
\begin{align*}
\E\|T_n\| &\ge 2^{-1}\|\sum_{i=0}^{n-1} A_i^{(n)}
X_i\varepsilon_i\| \ge 2^{-1}\E|X_0|\|\sum_{i=0}^{n-1} A_i^{(n)}
\varepsilon_i\| \\
&\ge \tilde{c}\E|X_0|\sqrt{n\log n} = c\sqrt{n\log n},
\end{align*}
where $\tilde{c}$ is universal and $c= \tilde{c}\E|X_0|$.

Let $T_n'$ be the Toeplitz matrix corresponding to
$X_0',\ldots,X_{n-1}'$. Since $X_n = X_n'$ for large $n$ and by
the above estimates $\E\|T_n\| \to \infty$, it is enough to prove
that with probability 1,
\begin{displaymath}
\lim_n \frac{\|T_n'\|}{\E\|T_n\|} = 1.
\end{displaymath}

We also have
\begin{align}\label{aux6}
\Big|\E\|T_n\| - \E\|T_n'\|\Big| & \le \E\Big|\|T_n\| -
\|T_n'\|\Big| \le \E \|T_n - T_n'\|\nonumber\\
& = \E \| \sum_{i=0}^{n-1} A_i^{(n)} X_i\ind{|X_i| > \sqrt{i\log i}}\|\nonumber\\
&\le 2\E|X_0| + 2\sum_{i=1}^{n-1} \E|X_i|\ind{|X_i| > \sqrt{i\log
i}}\nonumber\\
&\le 4 \E|X_0| + 2\E|X_0|^2 \sum_{i=2}^{n-1}\frac{1}{\sqrt{i\log
i}} = o(\sqrt{n\log n}) = o(\E\|T_n\|),
\end{align}
which shows that $\E\|T_n\|/\E\|T_n'\| \to 1$, so it is enough to
prove that with probability 1,
\begin{align}\label{reduction}
\lim_n\frac{\|T_n'\|}{\E\|T_n'\|} = 1.
\end{align}

For simplicity let us denote $s_n = \E \|T_n'\|$. Notice that by
the above inequalities, we have
\begin{align}\label{twosides}
c\sqrt{n\log n} \le s_n \le C\sqrt{n\log n}
\end{align}
for some $0 < c < C <\infty$.
 Fix $\theta > 1$ (but close to 1) and let $m_n$ be the
smallest integer $k$ such that $s_k > \theta^n$. We have $m_n \to
\infty$.

Notice that we have
\begin{displaymath}
s_{n+1} - s_n \le \E \sqrt{\sum_{i=0}^n 2 (X_i')^2} \le \sqrt{2n\E
X_0^2} = o(s_n).
\end{displaymath}
Thus
\begin{displaymath}
1< \frac{s_{m_n}}{\theta^n} \le \frac{s_{m_n}}{s_{m_n-1}} \le 1 +
\alpha_n,
\end{displaymath}
where $\alpha_n\to 0$, so $s_{m_n}/\theta^n \to 1$. In particular
$s_{m_{n+1}}/s_{m_n} \to \theta$.

Having established these preliminary facts, we can proceed with
the proof of (\ref{reduction}). We will show separately that
$\limsup \|T_n'\|/s_n \le 1$ and $\liminf \|T_n'\|/s_n \ge 1$ a.s.

We want to prove that for every $\varepsilon > 0$,
\begin{displaymath}
\lim_{n\to \infty}\p\Big(\sup_{k > n}\frac{\|T_k'\|}{s_k} > 1 +
\varepsilon\Big) = 0
\end{displaymath}
and
\begin{displaymath}
\lim_{n\to \infty} \p\Big(\inf_{k > n}\frac{\|T_k'\|}{s_k} < 1 -
\varepsilon\Big) = 0,
\end{displaymath} which will follow if for some $\theta > 1$ we show that
\begin{displaymath}
\sum_{k}\p\Big(\max_{ m_k \le i < m_{k+1}} \frac{\|T_i'\|}{s_i}
> 1 + \varepsilon\Big) < \infty
\end{displaymath}
and \begin{displaymath}
\sum_{k}\p\Big(\min_{ m_k \le i < m_{k+1}}
\frac{\|T_i'\|}{s_i} < 1 - \varepsilon\Big) < \infty
\end{displaymath}

Similarly as in the proof of Proposition \ref{non_iid}, we can now
take advantage of the monotonicity of $\|T_n'\|$ (and $s_n$),
which allows us to replace the above conditions by
\begin{align}\label{toprove}
\sum_{k}\p\Big(\|T_{m_{k+1}}'\|
> (1 + \varepsilon)s_{m_k}\Big) < \infty
\end{align}
and
\begin{align}\label{toprove1}
\sum_{k}\p\Big(\|T_{m_k}'\| < (1 - \varepsilon)s_{m_{k+1}}\Big) <
\infty
\end{align}

Let us now notice that since $\E X_i = 0$, we have
\begin{align*}
\|\E T_n '\| &=  \|\sum_{i=0}^{n-1} A_i^{(n)} \E X_i\ind{|X_i| \le
\sqrt{i\log i}}\| = \|\sum_{i=0}^{n-1} A_i^{(n)} \E X_i\ind{|X_i|
>
\sqrt{i\log i}}\|\\
&\le \E \|\sum_{i=0}^{n-1} A_i^{(n)} X_i\ind{|X_i| > \sqrt{i\log
i}}\| = o(\E\|T_n\|) = o(s_n),
\end{align*}
where in the third equality we used (\ref{aux6}).

Since $s_{m_{k+1}}/s_{m_k} \to \theta$, for large $k$,
\begin{align*}
\|T_{m_k}'\| &\ge \|T_{m_k}' - \E T_{m_k}'\| - \varepsilon
s_{m_{k+1}}/4, \\
\|T_{m_{k+1}}'\| &\le \|T_{m_{k+1}}' - \E T_{m_{k+1}}'\| +
\varepsilon s_{m_k}/4,
\end{align*}
so (\ref{toprove}) and (\ref{toprove1}) will follow if we show
that
\begin{align}\label{toprove2}
\sum_{k}\p\Big(\|T_{m_{k+1}}' - \E T_{m_{k+1}}' \|
> (1 + 3\varepsilon/4)s_{m_k}\Big) < \infty
\end{align}
and
\begin{align}\label{toprove3}
\sum_{k}\p\Big(\|T_{m_k}' - \E T_{m_k}'\| < (1 -
3\varepsilon/4)s_{m_{k+1}}\Big) < \infty.
\end{align}
Finally for $\theta \simeq 1$ (depending on $\varepsilon$) and
large $k$ we can write
\begin{align*}
\E\|T_{m_{k+1}}' - \E T_{m_{k+1}}' \| &\le s_{m_{k+1}}(1 +
\varepsilon/4) \le s_{m_k}(1 + \varepsilon/2),\\
\E\|T_{m_k}' - \E T_{m_k}'\| &\ge s_{m_k} (1 - \varepsilon/4) \ge
s_{m_{k+1}}(1-\varepsilon/2),
\end{align*}
so the proof of (\ref{toprove2}) and (\ref{toprove3}) can be
reduced to showing that
\begin{align*}
\sum_{k}\p\Big(\|T_{m_{k+1}}' - \E T_{m_{k+1}}' \|
> (1 + \eta_1)\E \|T_{m_{k+1}}' - \E T_{m_{k+1}}' \|  + \varepsilon s_{m_k}/8\Big) < \infty
\end{align*}
and
\begin{align*}
\sum_{k}\p\Big(\|T_{m_k}' - \E T_{m_k}'\| < (1 -
\eta_2)\E\|T_{m_k}' - \E T_{m_k}'\| -  \varepsilon
s_{m_{k+1}}/8\Big) < \infty,
\end{align*}
where
\begin{align*}
1 + \eta_1 & = (1 + 5\varepsilon/8)(1 + \varepsilon/2)^{-1},\\
1 -\eta_2 & = (1 - 5\varepsilon/8)(1 -\varepsilon/2)^{-1}.
\end{align*}

Now we can use Theorem \ref{FukNagaevThm} with $p = 3$ to obtain
\begin{align}\label{ineq1}
&\p\Big(\|T_{m_{k+1}}' - \E T_{m_{k+1}}' \|
> (1 + \eta_1)\E \|T_{m_{k+1}}' - \E T_{m_{k+1}}' \|  + \varepsilon
s_{m_k}/8\Big) \nonumber\\
&\le \exp\Big(-\frac{\varepsilon^2 s_{m_k}^2}{K m_{k+1}}\Big) +
K\frac{m_{k+1}\E|X_0|^3\ind{|X_0| \le \sqrt{m_{k+1}\log
m_{k+1}}}}{s_{m_k}^3\varepsilon^3}.
\end{align}

Similarly
\begin{align}\label{ineq2}
&\p\Big(\|T_{m_k}' - \E T_{m_k}'\| < (1 - \eta_2)\|T_{m_k}' - \E
T_{m_k}'\| -  \varepsilon s_{m_{k+1}}/8\Big)\nonumber\\
&\le K\exp\Big(-\frac{\varepsilon^2s_{m_{k+1}}^2}{Km_k}\Big) +
K\frac{m_k\E|X_0|^3\ind{|X_0| \le \sqrt{m_k\log
m_k}}}{s_{m_{k+1}}^3\varepsilon^3}.
\end{align}

By (\ref{twosides}) and the definition of $s_k$ we now obtain, that
for large $k$, $\sqrt{m_k\log m_k} \simeq \theta^k$  and $s_{m_k}
\simeq \theta^k$, so we also have $m_k \simeq \theta^{2k}/k$
(where by $\simeq$ we mean that the quotient of the left and right
hand side is bounded away from 0 and infinity by some constants,
perhaps depending on $\theta$, but independent of $k$). Thus the
right hand sides of (\ref{ineq1}) and (\ref{ineq2}) are bounded by
\begin{align*}
&K\exp\Big(-\frac{\varepsilon^2\theta^{2k}}{K\theta^{2k+2}k^{-1}}\Big)
+ K\frac{\theta^{2k+2}\E|X_0|^3\ind{|X_0| \le
C\theta^k}}{\theta^{3k}k\varepsilon^3} \\
&\le K\exp\Big(-\frac{\varepsilon^2 k}{K\theta^2}\Big) +
K\frac{\theta^{2}\E|X_0|^3\ind{|X_0| \le
C\theta^k}}{\theta^{k}k\varepsilon^3}.
\end{align*}
The series corresponding to the first term on the right hand side
is clearly convergent. As for the other term, we have
\begin{align*}
\sum_k \frac{\E|X_0|^3\ind{|X_0|\le C\theta^k}}{\theta^k k} &\le
\E|X_0|^3\sum_{k \colon \theta^k \ge |X_0|/C}  \theta^{-k}\le
\E|X_0|^3 \cdot \frac{KC}{|X_0|} \\
&= KC \E|X_0|^2 < \infty,
\end{align*}
which shows (\ref{toprove}) and (\ref{toprove1}) and thus proves
the theorem.
\end{proof}

\paragraph{Remark} Actually, in the above proof one may as well
use the inequality for bounded variables (e.g. Corollary
\ref{CLT_type}). It is enough to notice that all the other
calculations will still work if we truncate the variables $X_i$ at
the level $\sqrt{n}$ instead of $\sqrt{n\log n}$.

\paragraph{}As a corollary we obtain a law of large numbers for the i.i.d. case, proven originally in \cite{BS}.

\begin{cor}\label{lln_bs_iid}
Assume that $X_i$ are i.i.d. random variables, such that $\E X_0^2
< \infty$ and $\E X_0 = m$. Then
\begin{displaymath}
\lim_{n\to \infty} \frac{\|T_n\|}{n} = |m|
\end{displaymath}
\end{cor}

\begin{proof}
The argument we present comes from \cite{M}. We repeat it for the
sake of completeness. Let $T_n'$ be the matrix generated by the
variables $X_i -m$, $i= 0,\ldots,n-1$. By Theorems
\ref{expectation_from_above} and \ref{iid_case}, we have
\begin{displaymath}
\limsup_{n\to \infty} \frac{\|T_n'\|}{\sqrt{n\log n}} < \infty.
\end{displaymath}
But $T_n' = T_n - mI_n$, where $I_n$ is the $n\times n$ matrix of
ones. Since $\|I_n\| = n$, we have
\begin{align*}
\frac{\Big|\|T_n\|  - |m| n\Big|}{n} \le \frac{\|T_n'\|}{n} \to 0,
\end{align*}
which proves the corollary.
\end{proof}

\paragraph{Remark} A similar corollary for
Toeplitz matrices generated by independent but not identically
distributed random variables with the same mean may be derived
from Proposition \ref{non_iid}. Also, if we assume for instance
that $\E|X_n|$ is bounded away from $0$, we can combine the proofs
of Proposition \ref{non_iid} and Theorem \ref{iid_case} to obtain
the law of large numbers (with the normalization $\E\|T_n\|$) in
the centered, non i.i.d. case.

\paragraph{}Combining Theorem \ref{iid_case} with Lemma \ref{necessity}
gives the necessary and sufficient conditions for the $\limsup
\|T_n\|/\sqrt{n\log n}$ to be almost surely finite.

\begin{cor}\label{main_cor}
Assume that $X_i$ are i.i.d. random variables. Then
\begin{displaymath}
\limsup_{n \to \infty}\frac{\|T_n\|}{\sqrt{n\log n}} < \infty\;
{\rm a.s.},
\end{displaymath}
if and only if $\E X_0 = 0$ and $\E X_0^2 < \infty$.
\end{cor}

\begin{proof}
The first part follows from Theorems \ref{expectation_from_above}
and \ref{iid_case}. To prove the second part it is enough to show
that $\E X_0^2 < \infty$, since then Corollary \ref{lln_bs_iid}
gives $\E X_0 = 0$. To prove the finiteness of the second moment,
without loss of generality we may assume that the variable is
symmetric. The square integrability will follow by Lemma
\ref{necessity} if we prove that
\begin{displaymath}
\limsup_{n\to \infty}\frac{\E\|T_n\|}{\sqrt{n\log n}} < \infty.
\end{displaymath}
This can be proven in the same way as part (a) of Lemma 4.1 in
\cite{EL} (this lemma deals with partial sums of i.i.d. variables,
whereas in our case $T_n = \sum_{i=0}^{n-1} A_i^{(n)} X_i$,
however the fact that $\|A_i^{(n)}\|$ are bounded, uniformly in
$i$ and $n$, makes the same argument valid also in our situation).
\end{proof}

\paragraph{Problem} In view of results for other types of random
matrices, a natural question arises, whether
$\frac{\|T_n\|}{\sqrt{n\log n}}$ is almost surely convergent. From
Theorem \ref{iid_case} it follows that this is equivalent to the
existence of $\lim \frac{\E\|T_n\|}{\sqrt{n\log n}}$.

\subsection{Some concentration inequalities\label{last_section}}
We would like to conclude by presenting some concentration of
measure results for $\|T_n\|$.

Again, we will look at the matrix $T_n$ as a sum of independent
random variables in a Banach space. However, since the weak
variance $\sigma^2$, defined in Section \ref{concentration} is (in
the i.i.d. case) up to a universal constant the same as the so
called strong variance, defined in our situation as
\begin{align}\label{Sigma}
\Sigma^2 = \sum_{i=0}^{n-1}\E \|A_i\|^2 X_i^2 \le 4\sum_{i=1}^n \E
X_i^2,
\end{align}
we will pass to sums of real random variables. At the cost of
increasing some constants by a universal factor (in the i.i.d.
case), this will allow us to obtain concentration around
$\E\|T_n\|$ rather than bounds on the probability of deviation
from $(1\pm\eta)\E \|T_n\|$.

Our main technical tool will be the following symmetrization
lemma, proved in \cite{A} (it actually seems to be a part of the folklore,
but we haven't been able to find any references).
\begin{lemma}
\label{lemma_symmet} Let $\varphi \colon \R \to \R$ be a convex
function and $G = f(Y_1,\ldots,Y_n)$, where $Y_1,\ldots,Y_n$ are
independent random variables with values in a measurable space
$\mathcal{E}$. Denote
\begin{displaymath}
G_i = f(Y_1,\ldots,Y_{i-1},\tilde{Y}_i,Y_{i+1},\ldots,Y_n),
\end{displaymath}
where $(\tilde{Y}_1,\ldots, \tilde{Y}_n)$ is an independent copy
of $(Y_1,\ldots,Y_n)$ and assume that
\begin{displaymath}
|G - G_i| \le F_i(Y_i,\tilde{Y}_i)
\end{displaymath}
for some functions $F_i\colon \mathcal{E}^2 \to \R$, $i
=1,\ldots,n$. Then
\begin{displaymath}
\E \varphi(G - \E G) \le \E \varphi(\sum_{i=1}^n
\varepsilon_iF_i(Y_i,\tilde{Y}_i)), \label{symmet}
\end{displaymath}
where $\varepsilon_1,\ldots,\varepsilon_n$ is a sequence of
independent Rademacher variables, independent of $(Y_i)_{i=1}^n$
and $(\tilde{Y}_i)_{i=1}^n$.
\end{lemma}

Thus, whenever we have a concentration inequality for sums of
centered real variables, which has been obtained by (or can be
rephrased as) an estimate for expectations of one or many convex
functions (like the Lagrange transform or moment estimates), by
using the above Lemma we can obtain a corresponding result for
more general variables. In most cases the inequalities for sums of
Banach space valued variables obtained with this method will be
much worse than the optimal ones, since they lose the information
about the weak variance, in our case however the weak variance
satisfies
\begin{align*}
\sigma^2 &= \sup_{\|\alpha\|_2,\|\beta\|_2 \le 1} \sum_{k=0}^{n-1}
\Big(\sum_{i,j\colon\; |i-j| = k} \alpha_i\beta_j\Big)^2 \E X_k^2
\\&= \sup_{\|\alpha\|_2,\|\beta\|_2,\|\gamma_2\|_2 \le 1}
\Big( \sum_{k}  \sum_{i,j\colon\; |i-j| = k} \alpha_i\beta_j
\gamma_k (\E X_k^2)^{1/2}\Big)^2\\
&= \sup_{\|\gamma\|_2 \le 1} \Big\|[\gamma_{|i-j|} (\E
X^2_{|i-j|})^{1/2}]_{i,j\le n}\Big\|^2
\end{align*}
and in the i.i.d. case the choice of $\gamma_k = 1/\sqrt{n}$ for
all $k$ gives $\sigma^2 \ge n\E X_0^2$.

In the following proposition we present inequalities for Toeplitz
matrices generated by various classes of independent variables.
The proof relies on translating proper inequalities for real
variables into a moment form (which in this case can be easily
achieved with integration by parts), using Lemma
\ref{lemma_symmet} and then going back to tail estimates by
optimizing in $p$ the Chebyshev inequality for $p$-th moments. We
will skip the details which are quite standard but lengthy. The
real line inequalities we use are the $\psi_2$ inequality for
centered variables and Theorem \ref{unbounded}.

\begin{prop}\label{conc_ineq} Let $X_0,\ldots,X_{n-1}$ be independent random
variables and let $T_n$ the corresponding random Toeplitz matrix.
Then
\begin{itemize}
\item[(a)] if $\|X_i\|_{\psi_2} < \infty$ for all $i$, then for
all $t \ge 0$,
\begin{displaymath}
\p\Big(\Big|\|T_n\| - \E\|T_n\|\Big|\ge t \Big) \le
K\exp\Big(-\frac{t^2}{K \sum_{i=0}^{n-1} \|X_i\|_{\psi_2}^2}\Big),
\end{displaymath}
where $K$ is a universal constant,
\item[(b)] if for some $\alpha \in (0,1]$, $\|X_i\|_{\psi_\alpha}
< \infty$ for all $i$, then
\begin{displaymath}
\p\Big(\Big|\|T_n\| - \E\|T_n\|\Big|\ge t \Big) \le
2\exp\Big[-\frac{1}{K_\alpha}\min\Big(\frac{t^2}{\Sigma^2},\frac{t}{\|\max_i|X_i|\|_{\psi_\alpha}}\Big)\Big],
\end{displaymath}
where $\Sigma^2$ is defined by formula (\ref{Sigma}) and
$K_\alpha$ depends only on $\alpha$.
\end{itemize}
\end{prop}

\paragraph{Remark}
Obviously one can apply the same scheme to other inequalities. We
would also like to remark that in the non i.i.d. case one may
prefer to use Theorem \ref{unbounded} for sums of independent
Banach space valued variables. At the moment we do not know
whether also in this case $\sigma^2$ and $\Sigma^2$ are of the
same order.

\subsection{Final remarks}

In the article we restrict our attention to symmetric Toeplitz
matrices. However, the same methods provide analogous results for
nonsymmetric Toeplitz matrices $[X_{i-j}]_{i\le i,j\le n}$ where
$X_i$, $i \in \Z$ are independent random variables. Therefore, the
methods apply also to random Hankel matrices $[X_{j+k-1}]_{1\le
j,k\le n}$, since (as noticed in \cite{BDJ}, Remark 1.2.), the
singular values of $H_n$ are the same as of the nonsymmetric
Toeplitz matrix, obtained from $H_n$ by reversing the order of the
rows. The results about the expected operator norm of the
nonsymmetric Toeplitz matrix translate thus directly. As for
Theorem \ref{iid_case}, one has to be a little bit more careful.
This result does not translate formally, since the Toeplitz
matrices corresponding to $H_n$ are generated by a permutation of
the original sequence of random variables. One can however still
use all the concentration of measure results and simply repeat the
proof.


\begin{thebibliography}{5}
\bibitem{AKL} de Acosta A., Kuelbs J., Ledoux, M.
An inequality for the law of the iterated logarithm. Probability
in Banach spaces, IV (Oberwolfach, 1982), 1--29, Lecture Notes in
Math., 990, Springer, Berlin, 1983. MR0707506.
\bibitem{A} Adamczak R. A tail inequality for suprema of unbounded empirical processes with applications to Markov
chains. Submitted.

\bibitem{B} Bai Z.D. Methodologies in spectral analysis of
large-dimensional random matrices, a review. Statist. Sinica, 9
(3), 611-677, 1999. MR1711663

\bibitem{BS} Bose A., Sen A. Spectral norm of random large
dimensional noncentral Toeplitz and Hankel matrices. Electron.
Comm. Probab., 12, 29-35, 2007. MR2284045

\bibitem{BDJ} Bryc W., Dembo A., Jiang T., Spectral measure of
large Hankel, Markov and Toeplitz matrices. Ann. Probab., 34 (1),
1-38, 2006. MR2206341

\bibitem{EL} Einmahl U., Li D., Characterization of LIL behavior
in Banach space. To appear in Trans. Amer. Math. Soc.

\bibitem{FN} Fuk, D. H., Nagaev, S. V.
Probabilistic inequalities for sums of independent random
variables. Teor. Verojatnost. i Primenen. 16, 660--675, 1971. MR0293695

\bibitem{GLZ} Gin\'{e} E., Lata\l a R., Zinn J., Exponential and
moment inequalities for $U$-statistics. In \textit{High
Dimensional Probability II}, 13-38. Progr. Probab. 47. Birkhauser,
Boston, MA, 2000. MR1857312.

\bibitem{HM} Hammond C., Miller S.J., Distribution of eigenvalues
for the ensemble of real symmetric Toeplitz matrices. J. Theoret.
Probab., 18, 537-566, 2005. MR2167641.

\bibitem{LT} Ledoux M., Talagrand M., Probability in Banach
Spaces.  Isoperimetry and processes., Volume 23 of Ergebnisse der
Mathematik und ihrer Grenzgebiete (3). Springer-Verlag, Berlin
1991. MR1102015

\bibitem{KR} Klein T., Rio E., Concentration around the mean for
maxima of empirical processes. Ann. Probab. 33, 1060-1077, 2005. MR2135312

\bibitem{M} Meckes M., On the spectral norm of a random Toeplitz
matrix, Elect. Comm. in Probab., 12, 315-325, 2007. MR2342710

\bibitem{T} \textsc{Talagrand M.} New concentration inequalities in product spaces.
Invent. Math. 126, 1996, no. 3, 505--563. MR1419006.

\bibitem{VW} \textsc{van der Vaart, Aad W., Wellner, Jon A.}
Weak convergence and empirical processes. With applications to
statistics. Springer Series in Statistics. Springer-Verlag, New
York, 1996. MR1385671
\end{thebibliography}
\end{document}